# Fuzzy Calculus with Noval Approach Using Fuzzy Functions

Purnima Pandit
*Department of Applied Mathematics,*
*Faculty of Technology & Engineering,*
*M. S. University of Baroda, Vadodara, Gujarat, India*
pkpandit@yahoo.com

Payal Singh *
*Department of Applied Sciences and Humanities,*
*Faculty of Engineering and Technology,*
*Parul University, Vadodara, Gujarat, India*
singhpayalmath@gmail.com

## Abstract

This article deals with complexity in fuzzy functions where input, output and mapping are also fuzzy. The study of complexity in fuzzy functions, especially when applied to fuzzy derivatives and fuzzy calculus, is likely to involve considerations of how the uncertainty in inputs and function mapping affects the behaviour of the derivative, integration, and other calculus operations. This can have significant implications for decision-making, modelling, and analysis in various fields.

To establish fuzzy calculus with such complexity, we redefine Modified Hukuhara derivative. Under this derivative, fuzzy Taylor's series and its convergence is proposed and proved. Lastly, we solve a fully fuzzy differential equation with initial condition without converting it into the crisp one.

## 1 Introduction

In Fuzzy calculus, we study the derivatives and integrals of a fuzzy function. Since the beginning of fuzzy theory, lots of work has been carried in the field of fuzzy calculus. Different types of fuzzy derivatives are proposed and used in solving real-life applications of engineering and science. Initially, D. Dubois and H. Prade [1], [2], [3] have explained fuzzy mappings and given the difference between different kinds of approaches of fuzzy mapping. They also discussed the fuzzy differentiation of fuzzy valued functions. Kaleva [4] has studied the properties of differentiability and integrability.

In 1983, Puri and Ralescu [5], gave the idea about differential of fuzzy function (Hukuhara derivative) which extends the differential of a set-valued function. In the beginning, authors [4], [6] solved differential equations with fuzzy initial conditions and discussed the uniqueness and existence of the FDE under the Hukuhara derivative. The drawback in using the Hukuhara derivative is that the fuzzy solution becomes unbounded with time. In [7], [34], strongly generalized-Hukuhara derivative and generalized-Hukuhara differentiability concepts are proposed and proved, which overcomes the drawback of the Hukuhara derivative. But the limitation in using these derivatives is that the uniqueness of the solution is lost. Purnima and Payal in [8], recently introduced Modified Hukuhara derivative (mH-derivative), which gives a unique and bounded fuzzy solution for all time.

* Corresponding author.
*E-mail addresses:* pkpandit@yahoo.com (Purnima Pandit), singhpayalmath@gmail.com (Payal Singh).

Under these above mentioned fuzzy derivatives, fuzzy dynamical systems are solved using different techniques like Numerical techniques [9], [10], [11], [12], [13], [14], [15], [16], [17], [18], [19] [20], Transformation techniques [21], [22], [23], [24], [25], Analytical techniques [26], [27], [28], [30], Neural Network technique [31], [29] and Semi-Analytical techniques [32]. Numerical technique converts the differential equations in system of linear equations, which may involve negative fuzzy numbers. System of linear fuzzy equations with negative fuzzy number is solved by Purnima [35].

In this present article, we give the visualization of the complexity involved in fuzzy derivatives wherein input is also fuzzy. Based on this result, we propose and prove fuzzy Taylor's series along with its convergence under Modified Hukuhara derivative. We solve a fuzzy differential equation by fuzzy Taylor's series.

The following section contains preliminary, next to it is the section with various lemmas and main theorem along with their proofs followed with the illustrative example.

## 2 Preliminary

### 2.1 Fuzzy Number

$\mathcal{F} = \{\tilde{u}: R \to [0\ 1]\}$ is said to be a collection of a fuzzy number, if $\tilde{u}$ satisfies following properties.

- $\tilde{u}$ is normal.
- $\tilde{u}$ is a fuzzy convex.
- $\tilde{u}$ is upper semicontinuous.
- $\overline{\text{supp}\ (\tilde{u})} = \overline{\{x \in R / u(x) \geq 0\}}$ is compact.

### 2.2 Fuzzy Number in parametric form

$\alpha - cut$ is an important tool to convert fuzzy number into the crisp one and it is defined for

$0 < \alpha \leq 1,$ as ${}^{\alpha}\tilde{u} = \{\, x \in R/u(x) \ \geq \ \alpha\}.$ In general, it is defined in terms of parameter $\alpha$.

A fuzzy number in parametric form, obtained by performing $\alpha - cut,$ is an ordered pair of the form ${}^{\alpha}\tilde{A} = [\underline{A}, \overline{A}]$, satisfying the following conditions:

- $\underline{A}$ is bounded left continuous increasing function in $[0\ 1]$.
- $\overline{A}$ is bounded right continuous decreasing function in $[0\ 1]$.
- $\underline{A} \leq \overline{A}$.

### 2.3 Triangular Fuzzy Number

The triangular fuzzy number is denoted as triplet $(d, e, f)$, and its membership function $\tilde{A}(x)$, is given as

$$\tilde{A}(x) = \begin{cases} \frac{x-d}{e-d}, d < x \le e \\ \frac{f-x}{f-e}, e < x \le f \\ 0 \quad\quad otherwise \end{cases}$$

and it's $\alpha - cut$, is given as, ${}^{\alpha}\tilde{A} = [\underline{A}, \overline{A}] = [d + (e-d)\alpha, f - (f-e)\alpha]$.

## 2.4 Fuzzy arithmetic operations

Let $\tilde{A}$ and $\tilde{B}$ be two fuzzy numbers and $\lambda$ be any scalar. Scalar multiplication of fuzzy number is given as,

- $\lambda\, {}^{\alpha}\tilde{A} = \lambda[\underline{A}, \overline{A}] = [\lambda\underline{A}, \lambda\overline{A}]$.

Also, the arithmetic operations between $\tilde{A}$ and $\tilde{B}$ is defined using their parametric form, as follows.

- ${}^{\alpha}\tilde{A} \oplus {}^{\alpha}\tilde{B} = [\underline{A}, \overline{A}] \oplus [\underline{B}, \overline{B}\,] = [\underline{A} + \underline{B}, \overline{A} + \overline{B}]$
- ${}^{\alpha}\tilde{A} \otimes {}^{\alpha}\tilde{B} = [\underline{A}, \overline{A}] \otimes [\underline{B}, \overline{B}\,] = [\min(\underline{A}\,\underline{B}, \underline{A}\,\overline{B}, \overline{A}\,\underline{B}, \overline{A}\,\overline{B}), \max(\underline{A}\,\underline{B}, \underline{A}\,\overline{B}, \overline{A}\,\underline{B}, \overline{A}\,\overline{B})]$
- $\frac{{}^{\alpha}\tilde{A}}{{}^{\alpha}\tilde{B}} = [\underline{A}, \overline{A}] \otimes \left[\frac{1}{\overline{B}}, \frac{1}{\underline{B}}\right]$

## 2.5 Hausdorff distance [22]

The Hausdorff distance between the two fuzzy numbers is defined by,

$$d: \mathcal{F} \times \mathcal{F} \to R_{+} \cup \{0\}$$

$$d({}^{\alpha}\tilde{u}, {}^{\alpha}\tilde{v}) = \sup_{\alpha \in [0,1]} \max\{|\underline{u} - \underline{v}|, |\overline{u} - \overline{v}|\}$$

$d$ is metric in $\mathcal{F}$ which satisfies the following properties,

$d({}^{\alpha}\tilde{u} \oplus {}^{\alpha}\tilde{w}, {}^{\alpha}\tilde{v} \oplus {}^{\alpha}\tilde{w}) = d({}^{\alpha}\tilde{u}, {}^{\alpha}\tilde{v}), \ \forall\, \tilde{u}, \tilde{v}, \tilde{w} \in \mathcal{F}.$

$d(k\, {}^{\alpha}\tilde{u}, k\, {}^{\alpha}\tilde{v}) = |k| d({}^{\alpha}\tilde{u}, {}^{\alpha}\tilde{v})$

$d({}^{\alpha}\tilde{u} \oplus {}^{\alpha}\tilde{v}, {}^{\alpha}\tilde{w} \oplus {}^{\alpha}\tilde{e}) \le d({}^{\alpha}\tilde{u}, {}^{\alpha}\tilde{v}) \oplus d({}^{\alpha}\tilde{w}, {}^{\alpha}\tilde{e}), \forall\, \tilde{u}, \tilde{v}, \tilde{w}, \tilde{e} \in \mathcal{F}.$

then $(\mathcal{F}, d)$ is complete metric space.

## 2.6 Generalized Hukuhara Derivative [5]

Generalized Hukuhara (gH) difference is defined as follows, considering $K(X)$ as space of a nonempty convex and compact set of $X$, taking $A, B \in K(X)$ then gH difference is,

$$A \ominus B \Leftrightarrow \begin{cases} A = B + C \\ Or\ B = A + (-1)C \end{cases}$$

Let ${}^{\alpha}A = [\underline{a}, \overline{a}]$, ${}^{\alpha}B = [\underline{b}, \overline{b}]$, ${}^{\alpha}\, C = [\underline{c}, \overline{c}]$,

$$[\underline{a},\overline{a}] \ominus [\underline{b},\overline{b}] = [\underline{c},\overline{c}] = \begin{cases} i)\underline{a} = \underline{b} + \underline{c} \\ \quad \overline{a} = \overline{b} + \overline{c} \\ OR_{ii)}\ \underline{b} = \underline{a} - \overline{c} \\ \quad \overline{b} = \overline{a} - \underline{c} \end{cases}$$

So that $[\underline{a},\overline{a}] \ominus [\underline{b},\overline{b}] = [\underline{c},\overline{c}]$ always defined by,

$$[\underline{c},\overline{c}] = [\min(\underline{a}-\underline{b},\overline{a}-\overline{b}), \max(\underline{a}-\underline{b},\overline{a}-\overline{b})]$$

$'\ominus'$ refers gH- difference is defined as above.

A function $\tilde{f}: (a,b) \to E$ is said to be generalized Hukuhara differentiable at $t_0 \in (a, b)$, if

$$\dot{\tilde{f}}(t_0) = \lim_{h\to 0+} \frac{\tilde{f}(t_0+h) \ominus \tilde{f}(t_0)}{h}$$

## 3 Fuzzy Calculus

For proposing fuzzy Taylor's series, various definitions and results on fuzzy function are required. These results like fuzzy valued function for fuzzy input, fuzzy continuity, fuzzy derivative, and fuzzy power series along with its convergence.

### 3.1 Fuzzy function

Consider a fuzzy valued scalar function with fuzzy argument $\tilde{f}: \mathcal{F} \to \mathcal{F}$ . Its parametric form can be defined as follows,

$$^{\alpha}\tilde{f}(\tilde{x}) = [\underline{f}(\tilde{x}), \overline{f}(\tilde{x})], \text{ where, } \underline{f}(\tilde{x}) = \min \tilde{f}(\tilde{x}) \text{ and } \overline{f}(\tilde{x}) = \max \tilde{f}(\tilde{x})$$

Further, $^{\alpha}\tilde{f}(\tilde{x}) = \left[\underline{\underline{f}}(\underline{x},\overline{x}), \overline{\overline{f}}(\underline{x},\overline{x})\right]$

$$\text{where, } \underline{\underline{f}}(\underline{x},\overline{x}) = \min\left(\underline{f}(\underline{x},\overline{x}), \overline{f}(\underline{x},\overline{x})\right), \ \overline{\overline{f}}(\underline{x},\overline{x}) = \max\left(\underline{f}(\underline{x},\overline{x}), \overline{f}(\underline{x},\overline{x})\right) \quad (1)$$

This definition can be extended to $n$ dimensional fuzzy valued function, by now considering $\tilde{f}$ as, $\tilde{f}: \mathcal{F}^n \to \mathcal{F}^n$ where, $\tilde{f} = \{\tilde{f}_1(\tilde{x}_1, \tilde{x}_2, \dots \tilde{x}_n), \tilde{f}_2(\tilde{x}_1, \tilde{x}_2, \dots \tilde{x}_n), \dots \tilde{f}_n(\tilde{x}_1, \tilde{x}_2, \dots \tilde{x}_n)\}$ .

where parametric form,

$$^{\alpha}\tilde{f}_i(\tilde{x}_j) = \left[\underline{f_i}(\tilde{x}_1, \tilde{x}_2, \dots, \tilde{x}_n), \overline{f}_i(\tilde{x}_1, \tilde{x}_2, \dots, \tilde{x}_n)\right] \ \forall\, i = 1,2,3,\dots,n.$$

Further,

$$^{\alpha}\tilde{f}_i(\tilde{x}_j) = \left[\underline{\underline{f_i}}(\underline{x}_1, \underline{x}_2, \dots \underline{x}_n, \overline{x}_1, \overline{x}_2, \dots \overline{x}_n), \overline{\overline{f}}_i(\underline{x}_1, \underline{x}_2, \dots \underline{x}_n, \overline{x}_1, \overline{x}_2, \dots \overline{x}_n)\right], \forall i = 1,2,3,..,n$$

where, $\underline{\underline{f_i}}$ and $\overline{\overline{f}}_i$ are similar to as defined in equation (1) can be written as,

$$\underline{\underline{f_i}} = \min\left(\underline{f_i}(\underline{x}_1, \underline{x}_2, \dots \underline{x}_n), \overline{f}_i(\underline{x}_1, \underline{x}_2, \dots \underline{x}_n)\right), \forall i = 1,2,3,..,n$$

$$\overline{\overline{f}}_i = \max\left(\underline{f_i}(\underline{x_1},\underline{x_2},\dots\underline{x_n}),\overline{f}_i(\underline{x_1},\underline{x_2},\dots\underline{x_n})\right), \forall i = 1,2,3,..,n.$$

## 3.2 Fuzzy continuity

Let $\tilde{f}:\mathcal{F} \to \mathcal{F}$ be a fuzzy valued function and $\tilde{x},\tilde{x}_0 \in \mathcal{F}$ then $\tilde{f}$ is fuzzy continuous at a point $\tilde{x}_0$, if for any fixed number $\alpha \in (0\ 1]$ and any $\epsilon > 0$, $\exists\ \delta(\epsilon,\alpha)$ such that $d(\tilde{f}(\tilde{x}),\tilde{f}(\tilde{x}_0)) < \epsilon$ and $d(\tilde{x},\tilde{x}_0) < \delta(\epsilon,\alpha)$.

## 3.3 Modified Hukuhara derivative

In [8], we defined modified Hukuhara derivative for $\tilde{f}: I \to \mathcal{F}$, where $I$ is time interval. Now we extend this definition for the function $\tilde{f}:\mathcal{F} \to \mathcal{F}$.

A function $\tilde{f}:\mathcal{F} \to \mathcal{F}$ is said to be modified Hukuhara differentiable for an element $\dot{\tilde{f}}(\tilde{x}_0) \in \mathcal{F}$, such that for small $h > 0$, $\tilde{f}(\tilde{x}_0 + h) \ominus \tilde{f}(\tilde{x}_0)$, $\tilde{f}(\tilde{x}_0) \ominus \tilde{f}(\tilde{x}_0 - h)$ should exist and

$$\lim_{h\to 0+}\frac{\tilde{f}(\tilde{x}_0 + h) \ominus \tilde{f}(\tilde{x}_0)}{h} = \lim_{h\to 0-}\frac{\tilde{f}(\tilde{x}_0) \ominus \tilde{f}(\tilde{x}_0 - h)}{h} = \dot{\tilde{f}}(\tilde{x}_0) \qquad (2)$$

The equivalent parametric form for the first limit is given as,

$$\lim_{h\to 0+}\frac{{}^{\alpha}\tilde{f}(\tilde{x}_0 + h) \ominus {}^{\alpha}\tilde{f}(\tilde{x}_0)}{h}$$

$$=\lim_{h\to 0}\frac{1}{h}\left[\min\left(\underline{f}(\tilde{x}_0 + h) - \underline{f}(\tilde{x}_0),\overline{f}(\tilde{x}_0 + h) - \overline{f}(\tilde{x}_0)\right),\max\left(\underline{f}(\tilde{x}_0 + h) - \underline{f}(\tilde{x}_0),\overline{f}(\tilde{x}_0 + h) - \overline{f}(\tilde{x}_0)\right)\right]$$

$$=\lim_{h\to 0}\frac{1}{h}\left[\min\left(\underline{\underline{f}}(\underline{x}_0 + h,\overline{x}_0 + h) - \underline{\underline{f}}(\underline{x}_0,\overline{x}_0),\overline{\overline{f}}(\underline{x}_0 + h,\overline{x}_0 + h) - \overline{\overline{f}}(\underline{x}_0,\overline{x}_0)\right),\max\left(\underline{\underline{f}}(\underline{x}_0 + h,\overline{x}_0 + h) - \underline{\underline{f}}(\underline{x}_0,\overline{x}_0),\overline{\overline{f}}(\underline{x}_0 + h,\overline{x}_0 + h) - \overline{\overline{f}}(\underline{x}_0,\overline{x}_0)\right)\right]$$

Similarly, the second limit in equation (2) can be given as,

$$\lim_{h\to 0+}\frac{{}^{\alpha}\tilde{f}(\tilde{x}_0) \ominus {}^{\alpha}\tilde{f}(\tilde{x}_0 - \mathrm{h})}{h}$$

$$= \lim_{h\to 0}\frac{1}{h}\left[\min\left(\underline{f}(\tilde{x}_0) - \underline{f}(\tilde{x}_0 - h),\overline{f}(\tilde{x}_0) - \overline{f}(\tilde{x}_0 - h)\right),\max\left(\underline{f}(\tilde{x}_0) - \underline{f}(\tilde{x}_0 - h),\overline{f}(\tilde{x}_0) - \overline{f}(\tilde{x}_0 - h)\right)\right]$$

$$= \lim_{h\to 0}\frac{1}{h}\left[\min\left(\underline{\underline{f}}(\underline{x}_0,\overline{x}_0) - \underline{\underline{f}}(\underline{x}_0 - h,\overline{x}_0 - h),\overline{\overline{f}}(\underline{x}_0,\overline{x}_0) - \overline{\overline{f}}(\underline{x}_0 - h,\overline{x}_0 - h)\right),\max\left(\underline{\underline{f}}(\underline{x}_0,\overline{x}_0) - \underline{\underline{f}}(\underline{x}_0 - h,\overline{x}_0 - h),\overline{\overline{f}}(\underline{x}_0,\overline{x}_0) - \overline{\overline{f}}(\underline{x}_0 - h,\overline{x}_0 - h)\right)\right]$$

In the following examples, we compute the derivative of a given function using the definition as in equation (2).

### 3.3.1 Example: Derivative of $\tilde{f}(\tilde{x}) = \tilde{x}^2$ at $\tilde{x}_0$, $\tilde{x}_0 > 0$.

From the definition of derivative as in equation (2), we get,

$$
\begin{aligned}
{}^{\alpha}\dot{\tilde{f}}(\tilde{x}_0) &= \lim_{h\to 0+} \frac{{}^{\alpha}\tilde{f}(\tilde{x}_0+h) - {}^{\alpha}\tilde{f}(\tilde{x}_0)}{h} \\
&= \lim_{h\to 0+} \frac{\left[\underline{x}_0+h, \overline{x}_0+h\right]^2 - \left[\underline{x}_0, \overline{x}_0\right]^2}{h} \\
&= \lim_{h\to 0+} \frac{\left[\underline{x}_0+h, \overline{x}_0+h\right]\left[\underline{x}_0+h, \overline{x}_0+h\right] - \left[\underline{x}_0, \overline{x}_0\right]\left[\underline{x}_0, \overline{x}_0\right]}{h} \\
&= \lim_{h\to 0+} \frac{[x_l, x_u] - [x_{l0}, x_{u0}]}{h} \qquad (3)
\end{aligned}
$$

where,

$$x_l = \min\left(\left(\underline{x}_0+h\right)^2, \left(\underline{x}_0+h\right)(\overline{x}_0+h), (\overline{x}_0+h)^2\right)$$

$$x_u = \max\left(\left(\underline{x}_0+h\right)^2, \left(\underline{x}_0+h\right)(\overline{x}_0+h), (\overline{x}_0+h)^2\right)$$

$$x_{l0} = \min\left(\left(\underline{x}_0\right)^2, \left(\underline{x}_0\overline{x}_0\right), (\overline{x}_0)^2\right)$$

$$x_{u0} = \max\left(\left(\underline{x}_0\right)^2, \left(\underline{x}_0\overline{x}_0\right), (\overline{x}_0)^2\right)$$

Since, $\tilde{x}_0 > 0, x_l = \left(\underline{x}_0+h\right)^2, x_u = (\overline{x}_0+h)^2, x_{l0} = \left(\underline{x}_0\right)^2, x_{u0} = (\overline{x}_0)^2.$ (4)

Using equations (4) and (3).

$${}^{\alpha}\dot{\tilde{f}}(\tilde{x}_0) = \lim_{h\to 0+} \frac{[\min(x_l - x_{l0}, x_u - x_{u0}), \max(x_l - x_{l0}, x_u - x_{u0})]}{h}$$

$${}^{\alpha}\dot{\tilde{f}}(\tilde{x}_0) = \lim_{h\to 0+} \frac{\left[\left(\underline{x}_0+h\right)^2 - \left(\underline{x}_0\right)^2, (\overline{x}_0+h)^2 - (\overline{x}_0)^2\right]}{h}$$

$${}^{\alpha}\dot{\tilde{f}}(\tilde{x}_0) = \lim_{h\to 0+} \frac{\left[2\underline{x}_0 h + h^2, 2\overline{x}_0 h + h^2\right]}{h}$$

Similarly, ${}^{\alpha}\dot{\tilde{f}}(\tilde{x}_0) = \lim_{h\to 0+} \frac{{}^{\alpha}\tilde{f}(\tilde{x}_0) \ominus {}^{\alpha}\tilde{f}(\tilde{x}_0 - \mathrm{h})}{h} = \left[2\underline{x}_0, 2\overline{x}_0\right]$

$$\therefore\ {}^{\alpha}\dot{\tilde{f}}(\tilde{x}_0) = \left[2\underline{x}_0, 2\overline{x}_0\right]$$

Thus, by the Decomposition theorem as in, Section 3.2, $\dot{\tilde{f}}(\tilde{x}_0) = 2\tilde{x}_0$.

**3.3.2 Example:** Derivative of $\tilde{\boldsymbol{f}}(\tilde{\boldsymbol{x}}) = \boldsymbol{a}\tilde{\boldsymbol{x}}^{\boldsymbol{n}}, \boldsymbol{a} > \boldsymbol{0}$ at $\tilde{\boldsymbol{x}}_{\boldsymbol{0}}$ , $\tilde{\boldsymbol{x}}_{\boldsymbol{0}} > \boldsymbol{0}$.

By using definition as in equation (2), we have

$$ {}^{\alpha}\dot{\tilde{f}}(\tilde{x}_0) = \lim_{h\to 0+} \frac{{}^{\alpha}\tilde{f}(\tilde{x}_0 + h) - {}^{\alpha}\tilde{f}(\tilde{x}_0)}{h} $$

$$ = \lim_{h\to 0+} \frac{a[\underline{x}_0 + h, \overline{x}_0 + h]^n - a[\underline{x}_0, \overline{x}_0]^n}{h} $$

where,

$$ x_l = \min\left\{a(\underline{x}_0 + h)^n, a(\overline{x}_0 + h)(\underline{x}_0 + h)^{n-1}, a(\overline{x}_0 + h)^2(\underline{x}_0 + h)^{n-2} \dots, a(\overline{x}_0 + h)^n\right\}, $$

$$ x_u = \max\left\{a(\underline{x}_0 + h)^n, a(\overline{x}_0 + h)(\underline{x}_0 + h)^{n-1}, a(\overline{x}_0 + h)^2(\underline{x}_0 + h)^{n-2} \dots, a(\overline{x}_0 + h)^n\right\}, $$

$$ x_{l0} = \min\{a(\underline{x}_0)^n, a(\underline{x}_0{}^{n-1}\overline{x}_0), \dots, a(\overline{x}_0)^n\}, $$

$$ x_{u0} = \max\{a(\underline{x}_0)^n, a(\underline{x}_0{}^{n-1}\overline{x}_0), \dots, a(\overline{x}_0)^n\}. $$

Since, $\tilde{x}_0 > 0, x_l = a(\underline{x}_0 + h)^n, x_u = a(\overline{x}_0 + h)^n, x_{l0} = a(\underline{x}_0)^n, x_{u0} = a(\overline{x}_0)^n.$ (5)

Putting equation (5) in equation (3).

$$ {}^{\alpha}\dot{\tilde{f}}(\tilde{x}_0) = \lim_{h\to 0+} \frac{[a(\underline{x}_0 + h)^n - a(\underline{x}_0)^n, a(\overline{x}_0 + h)^n - a(\overline{x}_0)^n]}{h} $$

Also,

$$ {}^{\alpha}\dot{\tilde{f}}(\tilde{x}_0) = \lim_{h\to 0-} \frac{[a(\underline{x}_0 - h)^n - a(\underline{x}_0)^n, a(\overline{x}_0 - h)^n - a(\overline{x}_0)^n]}{h} $$

$$ {}^{\alpha}\dot{\tilde{f}}(\tilde{x}_0) = [na\underline{x}_0{}^{n-1}, na\overline{x}_0{}^{n-1}] $$

Thus, by the Decomposition theorem as in, Section 3.2, $\dot{\tilde{f}}(\tilde{x}_0) = na\tilde{x}_0{}^{n-1}$.

Using this definition of a fuzzy derivative, we give fuzzy power series along with its convergence results and substantiate it with examples.

### 3.4 Fuzzy power series and its convergence

We know that power series around point $x_0$, in crisp form is given as,

$$ \sum_{n=0}^{\infty} a_n(x - x_0)^n = a_0 + a_1(x - x_0)^1 + a_2(x - x_0)^2 + \cdots $$

A power series of fuzzy valued functions around the point $\tilde{x}_0$ can be given as

$$ \sum_{n=0}^{\infty} \tilde{a}_n \otimes (\tilde{x} \ominus \tilde{x}_0)^n = \tilde{a}_0 \oplus \tilde{a}_1 \otimes (\tilde{x} \ominus \tilde{x}_0) \oplus \tilde{a}_2 \otimes (\tilde{x} \ominus \tilde{x}_0)^2 \oplus \dots \quad (6) $$

where, $\tilde{a}_n$ are any fuzzy coefficients and *n* is a positive integer. The parametric form for equation (6) can be given as,

$$\sum_{n=0}^{\infty} {}^{\alpha}\tilde{a}_n \otimes {}^{\alpha}(\tilde{x} \ominus \tilde{x}_0)^n$$

$$= [\underline{a}_0, \overline{a}_0] \oplus [\underline{a}_1, \overline{a}_1] \otimes \left[[\underline{x}, \overline{x}] - [\underline{x}_0, \overline{x}_0]\right] \oplus [\underline{a}_2, \overline{a}_2] \otimes \left[[\underline{x}, \overline{x}] - [\underline{x}_0, \overline{x}_0]\right]^2$$

$$\dots \oplus [\underline{a}_n, \overline{a}_n] \otimes \left[[\underline{x}, \overline{x}] - [\underline{x}_0, \overline{x}_0]\right]^n \oplus \dots$$

$$= [\underline{a}_0, \overline{a}_0] \oplus [\underline{a}_1, \overline{a}_1] \otimes [\underline{x} - \underline{x}_0, \overline{x} - \overline{x}_0] \oplus [\underline{a}_2, \overline{a}_2] \otimes [\underline{x} - \underline{x}_0, \overline{x} - \overline{x}_0]^2$$

$$\dots \oplus [\underline{a}_n, \overline{a}_n] \otimes [\underline{x} - \underline{x}_0, \overline{x} - \overline{x}_0]^n \oplus \dots$$

$$= [\underline{a}_0, \overline{a}_0] \oplus [\underline{a}_1, \overline{a}_1] \otimes [\underline{x} - \underline{x}_0, \overline{x} - \overline{x}_0]$$

$$\oplus [\underline{a}_2, \overline{a}_2] \otimes \left[\min\left\{(\underline{x} - \underline{x}_0)^2, (\overline{x} - \overline{x}_0)(\underline{x} - \underline{x}_0), (\overline{x} - \overline{x}_0)^2\right\}, \max\left\{(\underline{x} - \underline{x}_0)^2, (\overline{x} - \overline{x}_0)(\underline{x} - \underline{x}_0), (\overline{x} - \overline{x}_0)^2\right\}\right] \dots$$

$$\oplus [\underline{a}_n, \overline{a}_n] \otimes \left[\min\left\{\left((\underline{x} - \underline{x}_0)^n, (\overline{x} - \overline{x}_0)(\underline{x} - \underline{x}_0)^{n-1}, (\overline{x} - \overline{x}_0)^2(\underline{x} - \underline{x}_0),^{n-2} \dots (\overline{x} - \overline{x}_0)^n\right), (\overline{x} - \overline{x}_0)^n, (\overline{x} - \overline{x}_0)^{n-1}(\underline{x} - \underline{x}_0), (\overline{x} - \overline{x}_0)^{n-2}(\underline{x} - \underline{x}_0),^2 \dots, (\underline{x} - \underline{x}_0)^n\right\}, \max\left\{\left((\underline{x} - \underline{x}_0)^n, (\overline{x} - \overline{x}_0)(\underline{x} - \underline{x}_0)^{n-1}, (\overline{x} - \overline{x}_0)^2(\underline{x} - \underline{x}_0),^{n-2} \dots (\overline{x} - \overline{x}_0)^n\right), (\overline{x} - \overline{x}_0)^n, (\overline{x} - \overline{x}_0)^{n-1}(\underline{x} - \underline{x}_0), (\overline{x} - \overline{x}_0)^{n-2}(\underline{x} - \underline{x}_0),^2 \dots, (\underline{x} - \underline{x}_0)^n\right\}\right] \oplus ..$$

Since, $\tilde{x}_0 > 0$, then the above series expansion is given as follows,

$$\sum_{n=0}^{\infty} {}^{\alpha}\tilde{a}_n \otimes {}^{\alpha}(\tilde{x} \ominus \tilde{x}_0)^n$$

$$= [\underline{a}_0, \overline{a}_0] \oplus [\underline{a}_1, \overline{a}_1] \otimes [\underline{x} - \underline{x}_0, \overline{x} - \overline{x}_0] \oplus [\underline{a}_2, \overline{a}_2] \otimes \left[(\underline{x} - \underline{x}_0)^2, (\overline{x} - \overline{x}_0)^2\right] \dots$$

$$\oplus [\underline{a}_n, \overline{a}_n] \otimes \left[(\underline{x} - \underline{x}_0)^n, (\overline{x} - \overline{x}_0)^n\right] \oplus \dots$$

$$= [\underline{a}_0, \overline{a}_0]$$

$$\oplus \left[\min\{\underline{a}_1 (\underline{x} - \underline{x}_0), \underline{a}_1(\overline{x} - \overline{x}_0), \overline{a}_1(\underline{x} - \underline{x}_0), \overline{a}_1(\overline{x} - \overline{x}_0)\}, \max\{\underline{a}_1 (\underline{x} - \underline{x}_0), \underline{a}_1(\overline{x} - \overline{x}_0), \overline{a}_1(\underline{x} - \underline{x}_0), \overline{a}_1(\overline{x} - \overline{x}_0)\}\right]$$

$$\oplus \left[\min\left\{\underline{a}_1\left(\underline{x}-\underline{x}_0\right)^2, \underline{a}_1(\overline{x}-\overline{x}_0)^2, \overline{a}_1\left(\underline{x}-\underline{x}_0\right)^2, \overline{a}_1(\overline{x}-\overline{x}_0)^2\right\}, \max\left\{\underline{a}_1\left(\underline{x}-\underline{x}_0\right)^2, \underline{a}_1(\overline{x}-\overline{x}_0)^2, \overline{a}_1\left(\underline{x}-\underline{x}_0\right)^2, \overline{a}_1(\overline{x}-\overline{x}_0)^2\right\}\right]$$

$$\ldots$$

$$\oplus \left[\min\{\underline{a}_1\left(\underline{x}-\underline{x}_0\right)^n, \underline{a}_1(\overline{x}-\overline{x}_0)^n, \overline{a}_1\left(\underline{x}-\underline{x}_0\right)^n, \overline{a}_1(\overline{x}-\overline{x}_0)^n\}, \max\{\underline{a}_1\left(\underline{x}-\underline{x}_0\right)^n, \underline{a}_1(\overline{x}-\overline{x}_0)^n, \overline{a}_1\left(\underline{x}-\underline{x}_0\right)^n, \overline{a}_1(\overline{x}-\overline{x}_0)^n\}\right]$$

$$\oplus \ldots$$

$$= [\underline{a}_0, \overline{a}_0] \oplus [\underline{a}_1(\underline{x}-\underline{x}_0), \overline{a}_1(\overline{x}-\overline{x}_0)] \oplus \left[\underline{a}_2\left(\underline{x}-\underline{x}_0\right)^2, \overline{a}_2(\overline{x}-\overline{x}_0)^2\right] \ldots$$

$$\oplus \left[\underline{a}_n\left(\underline{x}-\underline{x}_0\right)^n, \overline{a}_n(\overline{x}-\overline{x}_0)^n\right] \oplus \ldots$$

In the next section, we give a result of a radius of convergence for fuzzy power series.

### 3.4.1 Radius of convergence

If $\widetilde{a_n} \neq 0$, a radius of convergence $\tilde{R}$ given by $\lim_{n\to\infty}\left|\frac{\widetilde{a_n}}{\widetilde{a_{n+1}}}\right|$ of fuzzy power series can be defined as in parametric form, ${}^{\alpha}\tilde{R} = [\underline{R}, \overline{R}]$, where,

$$\underline{R} = \min\left[\lim_{n\to\infty}\left|\frac{\underline{a}_n}{\underline{a}_{n+1}}\right|, \lim_{n\to\infty}\left|\frac{\underline{a}_n}{\overline{a}_{n+1}}\right|, \lim_{n\to\infty}\left|\frac{\overline{a}_n}{\underline{a}_{n+1}}\right|, \lim_{n\to\infty}\left|\frac{\overline{a}_n}{\overline{a}_{n+1}}\right|\right]$$

$$\overline{R} = \max\left[\lim_{n\to\infty}\left|\frac{\underline{a}_n}{\underline{a}_{n+1}}\right|, \lim_{n\to\infty}\left|\frac{\underline{a}_n}{\overline{a}_{n+1}}\right|, \lim_{n\to\infty}\left|\frac{\overline{a}_n}{\underline{a}_{n+1}}\right|, \lim_{n\to\infty}\left|\frac{\overline{a}_n}{\overline{a}_{n+1}}\right|\right]$$

The fuzzy number $\tilde{R}$ can be obtained by the Decomposition Theorem as given in Section **Error! Reference source not found.**.

The fuzzy power series $\sum \tilde{a}_n \otimes (\tilde{x} \ominus \tilde{x}_0)^n$ with a radius of convergence $\tilde{R}$ and the set of the points from an interval at which fuzzy series is convergent, known as the interval of convergence such that $d(\tilde{x}, \tilde{x}_0) < \tilde{R}$.

The parametric form of ${}^{\alpha}d(\tilde{x}, \tilde{x}_0) = \left[d\left(\underline{x}, \underline{x}_0\right), d(\overline{x}, \overline{x}_0)\right]$ and ${}^{\alpha}\tilde{R} = [\underline{R}, \overline{R}]$,

Then, $d(\tilde{x}, \tilde{x}_0) < \tilde{R}$ can be simplified in the following manner by writing in parametric form,

So, the component-wise expression is as follows,

$$d\left(\underline{x}, \underline{x}_0\right) < \underline{R}$$

$$-\underline{R} < \left(\underline{x} - \underline{x}_0\right) < \underline{R}$$

$$(\underline{x}_0 - \underline{R}) < \underline{x} < (\underline{x}_0 + \underline{R})$$

Similarly,

$$d(\overline{x}, \overline{x}_0) < \overline{R}$$

$$-\overline{R} < (\overline{x} - \overline{x}_0) < \overline{R}$$

$$(\overline{x}_0 - \overline{R}) < \overline{x} < (\overline{x}_0 + \overline{R})$$

Thus, when we combine the above result of a radius of convergence in parametric form, we obtain the following result for convergence,

$$\max d(\underline{x}, \underline{x}_0) < \underline{x} < \overline{x} < \min d(\overline{x}, \overline{x}_0)$$

This is the interval (sense of $\alpha - cut$) in which the fuzzy power series converges.

We prove the lemma for the convergence of power series based on the definition given in Section 3.4.1.

### 3.4.2 Lemma

If, $\underline{L} = \frac{1}{\underline{R}} = \min[\lim_{n\to\infty}\left|\frac{\underline{a}_{n+1}}{\underline{a}_n}\right|, \lim_{n\to\infty}\left|\frac{\underline{a}_{n+1}}{\overline{a}_n}\right|, \lim_{n\to\infty}\left|\frac{\overline{a}_{n+1}}{\underline{a}_n}\right|, \lim_{n\to\infty}\left|\frac{\overline{a}_{n+1}}{\overline{a}_n}\right|],$

$$\overline{L} = \frac{1}{\overline{R}} = \max[\lim_{n\to\infty}\left|\frac{\underline{a}_{n+1}}{\underline{a}_n}\right|, \lim_{n\to\infty}\left|\frac{\underline{a}_{n+1}}{\overline{a}_n}\right|, \lim_{n\to\infty}\left|\frac{\overline{a}_{n+1}}{\underline{a}_n}\right|, \lim_{n\to\infty}\left|\frac{\overline{a}_{n+1}}{\overline{a}_n}\right|]$$

and $\underline{L} < 1$ and $\overline{L} < 1$ then fuzzy power series in parametric form converges absolutely.

**Proof:** Given that,

$$\underline{L} = \min[\lim_{n\to\infty}\left|\frac{\underline{a}_{n+1}}{\underline{a}_n}\right|, \lim_{n\to\infty}\left|\frac{\underline{a}_{n+1}}{\overline{a}_n}\right|, \lim_{n\to\infty}\left|\frac{\overline{a}_{n+1}}{\underline{a}_n}\right|, \lim_{n\to\infty}\left|\frac{\overline{a}_{n+1}}{\overline{a}_n}\right|],$$

Let $\frac{b_{n+1}}{b_n}$ be the minimum value of $\min[\lim_{n\to\infty}\left|\frac{\underline{a}_{n+1}}{\underline{a}_n}\right|, \lim_{n\to\infty}\left|\frac{\underline{a}_{n+1}}{\overline{a}_n}\right|, \lim_{n\to\infty}\left|\frac{\overline{a}_{n+1}}{\underline{a}_n}\right|, \lim_{n\to\infty}\left|\frac{\overline{a}_{n+1}}{\overline{a}_n}\right|],$

and, $\left|\frac{b_{n+1}}{b_n}\right| < 1.$

Thus,

$$|b_{n+1}| < |rb_n|$$

$$|b_{n+i}| < |r^i b_n| \text{ for all } n > N \text{ and } i = 1,2,3,\ldots\ldots$$

Hence,

$$\sum_{i=N+1}^{\infty}|b_{n+i}| < \sum_{i=1}^{\infty}|r^i b_n|$$

The right-hand side of the above equation is convergent absolutely.

Similarly, series converges for $\overline{L} < 1$.

Hence, fuzzy power series in parametric form converges absolutely.

In the next section, we give an example for the expansion of a fuzzy power series along with its convergence.

### 3.4.3 Example:

The fuzzy power series for $\sum_{n=0}^{\infty} \widetilde{-1} \otimes \tilde{x}^n$ can be obtained as follows,

where, $\widetilde{-1} = (-2, -1, 1)$

$$\sum_{n=0}^{\infty} {}^{\alpha}\widetilde{-1} \otimes {}^{\alpha}\tilde{x}^n$$

$$= \sum_{n=1}^{\infty} [\alpha - 2, 1 - 2\alpha] \otimes [\underline{x^n}, \overline{x^n}]$$

$$= \sum_{n=0}^{\infty} \Big[\min\{(\alpha - 2) \otimes \underline{x^n}, (\alpha - 2) \otimes \overline{x^n}, (1 - 2\alpha) \otimes \underline{x^n}, (1 - 2\alpha) \otimes \overline{x^n}\}, \max\{(\alpha - 2) \otimes \underline{x^n}, (\alpha - 2) \otimes \overline{x^n}, (1 - 2\alpha) \otimes \underline{x^n}, (1 - 2\alpha) \otimes \overline{x^n}\}\Big] \tag{7}$$

Putting different values of *n* in equation (7), we can obtain the expression of the fuzzy power series. The first three terms in this expansion can be written as,

$$\sum_{n=0}^{\infty} {}^{\alpha}\widetilde{-1} \otimes {}^{\alpha}\tilde{x}^n$$

$$= [(\alpha - 2), (1 - 2\alpha)]$$

$$\oplus \Big[\min\{(\alpha - 2) \otimes \underline{x}^1, (\alpha - 2) \otimes \overline{x}^1, (1 - 2\alpha) \otimes \underline{x}^1, (1 - 2\alpha) \otimes \overline{x}^1\}, \max\{(\alpha - 2) \otimes \underline{x}^1, (\alpha - 2) \otimes \overline{x}^1, (1 - 2\alpha) \otimes \underline{x}^1, (1 - 2\alpha) \otimes \overline{x}^1\}\Big]$$

$$\oplus \Big[\min\{(\alpha - 2) \otimes \underline{x^2}, (\alpha - 2) \otimes \overline{x^2}, (\alpha - 2) \otimes \underline{x} \otimes \overline{x}, (1 - 2\alpha) \otimes \underline{x^2}, (1 - 2\alpha) \otimes \overline{x^2}, (1 - 2\alpha) \otimes \underline{x} \otimes \overline{x}\}, \max\{(\alpha - 2) \otimes \underline{x^2}, (\alpha - 2) \otimes \underline{x} \otimes \overline{x}\ (\alpha - 2) \otimes \overline{x^2}, (1 - 2\alpha) \otimes \underline{x^2}, (1 - 2\alpha) \otimes \overline{x^2}, (1 - 2\alpha) \otimes \underline{x} \otimes \overline{x}\}\Big]$$

$$\oplus \ldots \tag{8}$$

The radius of convergence for $\sum_{n=0}^{\infty} \widetilde{-1} \otimes \tilde{x}^n$ is defined as,

$$\tilde{R} = \lim_{n \to \infty} \left|\frac{\tilde{a}_n}{\tilde{a}_{n+1}}\right| = \lim_{n \to \infty} \left|\frac{\tilde{1}}{\tilde{1}}\right| = \tilde{1}.$$

And, the region of convergence by Lemma 3.4.2 is $|\tilde{x}| < \tilde{1}$.

The following section contains example of radius of convergence for fuzzy power series.

### 3.4.4 Example:

Find the radius of convergence of $\sum_{n=0}^{\infty} \tilde{2} \otimes \tilde{x}^n$, where $\tilde{2} = (1,2,3)$.

$$\tilde{R} = \lim_{n \to \infty} \left|\frac{\tilde{a}_n}{\tilde{a}_{n+1}}\right| = \lim_{n \to \infty} \frac{\tilde{2}}{\tilde{2}}$$

Let, ${}^{\alpha}\tilde{2} = [1 + \alpha, 3 - \alpha]$ then

$$\underline{R} = \min[\lim_{n\to\infty} \frac{\underline{a}_n}{\underline{a}_{n+1}}, \lim_{n\to\infty} \frac{\underline{a}_n}{\overline{a}_{n+1}}, \lim_{n\to\infty} \frac{\overline{a}_n}{\underline{a}_{n+1}}, \lim_{n\to\infty} \frac{\overline{a}_n}{\overline{a}_{n+1}}]$$

$$=\min\left[1, \frac{(1+\alpha)}{(3-\alpha)}, \frac{(3-\alpha)}{(1+\alpha)}, 1\right] = \frac{(1+\alpha)}{(3-\alpha)}$$

$$\overline{R} = \max[\lim_{n\to\infty} \frac{\underline{a}_n}{\underline{a}_{n+1}}, \lim_{n\to\infty} \frac{\underline{a}_n}{\overline{a}_{n+1}}, \lim_{n\to\infty} \frac{\overline{a}_n}{\underline{a}_{n+1}}, \lim_{n\to\infty} \frac{\overline{a}_n}{\overline{a}_{n+1}}]$$

$$= \max\left[1, \frac{(1+\alpha)}{(3-\alpha)}, \frac{(3-\alpha)}{(1+\alpha)}, 1\right] = \frac{(3-\alpha)}{(1+\alpha)}$$

Therefore, fuzzy radius of convergence is $\tilde{1}$, and in parametric form $\left[\frac{(1+\alpha)}{(3-\alpha)}, \frac{(3-\alpha)}{(1+\alpha)}\right]$. When we put $\alpha = 1,$ the radius of convergence is 1 which is same as its counter crisp problem.

The given fuzzy series is convergent for $|\tilde{x}| < \tilde{1}$, where $\tilde{1} = \left(\frac{1}{3}, 1, 3\right)$.

### 3.4.5 Example:

Consider the fuzzy power series $\sum \frac{\boldsymbol{n}}{\tilde{\mathbf{5}}^{\boldsymbol{n-1}}} \otimes (\tilde{\boldsymbol{x}} \oplus \tilde{\mathbf{2}})^{\boldsymbol{n}}$ where, $\tilde{\mathbf{2}} = (\mathbf{1}, \mathbf{2}, \mathbf{3})$ and $\tilde{\mathbf{5}} = (\mathbf{4}, \mathbf{5}, \mathbf{6})$.

Here, $\tilde{a}_n = \frac{n}{\tilde{5}^{n-1}}$.

The radius of convergence is given by,

$$\tilde{R} = \lim_{n\to\infty} \left|\frac{\tilde{a}_n}{\tilde{a}_{n+1}}\right| = \lim_{n\to\infty} \left|\frac{n}{\tilde{5}^{n-1}} \frac{\tilde{5}^n}{(n+1)}\right| = \tilde{5}$$

Thus, series is convergent for those $\tilde{x}$ which satisfies this condition $|\tilde{x} \oplus \tilde{2}| < \tilde{5}$

That is, $-\tilde{7} < \tilde{x} < \tilde{3}$.

The convergence of $|\tilde{x} \oplus \tilde{2}| < \tilde{5}$ can be visualized in crisp form by using the parametric form.

Let ${}^{\alpha}\tilde{x} = [\underline{x}\ \overline{x}]$, ${}^{\alpha}\tilde{2} = [1+\alpha,\ 3-\alpha]$ and ${}^{\alpha}\tilde{5} = [4+\alpha, 6-\alpha]$, then the inequality becomes

$$|[\underline{x}\ \overline{x}] + [1+\alpha, 3-\alpha]| < |4+\alpha, 6-\alpha|$$

That is,

$$|\underline{x} + 1 + \alpha| < 4 + \alpha$$

which implies,

$$-(4+\alpha) < (\underline{x} + 1 + \alpha) < 4 + \alpha$$

Thus,

$$-(5+2\alpha) < \underline{x} < 3$$

Similarly, $|\overline{x} + 3 - \alpha| < 6 - \alpha$

$$\Rightarrow -(6-\alpha) < (\overline{x} + 3 - \alpha) < 6 - \alpha$$

$$\Rightarrow (-9+2\alpha) < \overline{x} < 3$$

By 3.4.1, we can combine the results in crisp form,

$$[(-9+2\alpha), -(5+2\alpha)] < [\underline{x}, \overline{x}] < [3,3]$$

And by the Decomposition theorem given in section 3.2, we can write that the given series converges for the interval $-\tilde{7} < \tilde{x} < \tilde{3}$.
The above-mentioned theories, fuzzy power series and its convergence are required for proving Fuzzy Taylor's series.

In the next section, Fuzzy Taylor's series is proposed and proved.

### 3.5 Fuzzy Taylor's Theorem

If a fuzzy valued function, $\tilde{f}: \mathcal{F}^n \to \mathcal{F}^n$ is, *n* times modified Hukuhara (mh) differentiable, then fuzzy Taylor's expansion is given as,

$$\tilde{f}(\tilde{x}) = \tilde{f}(\tilde{x}_0) \oplus \dot{\tilde{f}}(\tilde{x}_0) \otimes (\tilde{x} \ominus \tilde{x}_0) \oplus \frac{\ddot{\tilde{f}}(\tilde{x}_0)}{2!} \otimes (\tilde{x} \ominus \tilde{x}_0)^2 \oplus \frac{\dddot{\tilde{f}}(\tilde{x}_0)}{3!} \otimes (\tilde{x} \ominus \tilde{x}_0)^3 \oplus \frac{\tilde{f}^4(\tilde{x}_0)}{4!} \otimes (\tilde{x} \ominus \tilde{x}_0)^4 \oplus \frac{\tilde{f}^5(\tilde{x}_0)}{5!} \otimes (\tilde{x} \ominus \tilde{x}_0)^5 \oplus \ldots\ldots\ldots \oplus \frac{\tilde{f}^n(\tilde{x}_0)}{n!} \otimes (\tilde{x} \ominus \tilde{x}_0)^n \oplus \ldots$$

**Proof:**

Let a fuzzy valued function $\tilde{f}$ is Modified Hukuhara differentiable then it can be approximated near a point $\tilde{x}_0$ by its tangent line that gives its linear approximation to $\tilde{f}$ at the point $\tilde{x}_0$.

The approximation is given by using the Modified Hukuhara derivative as in Section 3.5.

$$\dot{\tilde{f}}(\tilde{x}_0) = \lim_{\tilde{x} \to \tilde{x}_0} \frac{\tilde{f}(\tilde{x}) \ominus \tilde{f}(\tilde{x}_0)}{(\tilde{x} \ominus \tilde{x}_0)}$$

Now, in the neighbourhood $\tilde{x}_0$, discretization of first mh-derivative, gives us

$$\tilde{f}(\tilde{x}) \ominus \tilde{f}(\tilde{x}_0) = \dot{\tilde{f}}(\tilde{x}_0) \otimes (\tilde{x} \ominus \tilde{x}_0)$$

$$\therefore \tilde{f}(\tilde{x}) = \tilde{f}(\tilde{x}_0) \oplus \dot{\tilde{f}}(\tilde{x}_0) \otimes (\tilde{x} \ominus \tilde{x}_0)$$

If first approximation, $\tilde{p}_1(x)$ is polynomial of first degree and this polynomial has some properties, $\tilde{p}_1(\tilde{x}_0) = \tilde{f}(\tilde{x}_0), \dot{\tilde{p}}_1(\tilde{x}_0) = \tilde{f}'(\tilde{x}_0)$, then

$$\tilde{p}_1(\tilde{x}) = \tilde{f}(\tilde{x}_0) \oplus \dot{\tilde{f}}(\tilde{x}_0) \otimes (\tilde{x} \ominus \tilde{x}_0).$$

Linear approximation of $\tilde{f}$ is well defined at $\tilde{x}_0$, if $\tilde{f}$ has a constant slope. However, if $\tilde{f}$ has curvature near pt. $\tilde{x}_0$ then, it requires quadratic approximation. For quadratic approximation, we add one more term,

$$\tilde{p}_2(x) = \tilde{f}(\tilde{x}_0) \oplus \dot{\tilde{f}}(\tilde{x}_0) \otimes (\tilde{x} \ominus \tilde{x}_0) \oplus \tilde{c} \otimes (\tilde{x} \ominus \tilde{x}_0) \otimes (\tilde{x} \ominus \tilde{x}_0) \qquad (9)$$

To determine a new term $\tilde{c}$, $\tilde{p}_2(x)$ must be a good approximation to $\tilde{f}$, near the point $\tilde{x}_0$. This requires, $\tilde{p}_2(\tilde{x}_0) = \tilde{f}(\tilde{x}_0), \dot{\tilde{p}}_2(\tilde{x}_0) = \dot{\tilde{f}}(\tilde{x}_0)$ and $\ddot{\tilde{p}}_2(\tilde{x}_0) = \ddot{\tilde{f}}(\tilde{x}_0)$.

Differentiating, equation (9), and we get,

$$\dot{\tilde{p}}_2(\tilde{x}) = \dot{\tilde{f}}(\tilde{x}_0) \oplus 2\, \tilde{c} \otimes (\tilde{x} \ominus \tilde{x}_0)$$

Again differentiating, we get

$$\ddot{\tilde{p}}_2(\tilde{x}) = 2\,\tilde{c} \;\Rightarrow \tilde{c} = \frac{\ddot{\tilde{p}}_2(x)}{2}$$

At $\tilde{x} = \tilde{x}_0,\ \ddot{\tilde{p}}_2(\tilde{x}_0) = \ddot{\tilde{f}}(\tilde{x}_0)$

So, equation (9) becomes,

$$\tilde{p}_2(\tilde{x}) \ = \tilde{f}(\tilde{x}_0) \oplus \dot{\tilde{f}}(\tilde{x}_0) \otimes (\tilde{x} \ominus \tilde{x}_0) \oplus \frac{\ddot{\tilde{f}}(\tilde{x}_0)}{2} \otimes (\tilde{x} \ominus \tilde{x}_0) \otimes (\tilde{x} \ominus \tilde{x}_0)$$

Following a similar way, we get the third approximation as

$$\tilde{p}_3(\tilde{x}) = \tilde{f}(\tilde{x}_0) \oplus \dot{\tilde{f}}(\tilde{x}_0) \ \otimes (\tilde{x} \ominus \tilde{x}_0) \oplus \frac{\ddot{\tilde{f}}(\tilde{x}_0)}{2} \otimes (\tilde{x} \ominus \tilde{x}_0) \otimes (\tilde{x} \ominus \tilde{x}_0) \oplus \frac{\dddot{\tilde{f}}(\tilde{x}_0)}{6} \otimes (\tilde{x} \ominus \tilde{x}_0) \otimes (\tilde{x} \ominus \tilde{x}_0) \otimes (\tilde{x} \ominus \tilde{x}_0)$$

And the $n^{th}$ degree approximation to $\tilde{f}$, can be given as

$$\tilde{p}_n(\tilde{x}) = \tilde{f}(\tilde{x}_0) \oplus \dot{\tilde{f}}(\tilde{x}_0) \otimes (\tilde{x} \ominus \tilde{x}_0) \oplus \frac{\ddot{\tilde{f}}(\tilde{x}_0)}{2!} \otimes (\tilde{x} \ominus \tilde{x}_0)^2 \oplus \frac{\dddot{\tilde{f}}(\tilde{x}_0)}{3!} \otimes (\tilde{x} \ominus \tilde{x}_0)^3 \oplus \frac{\tilde{f}^4(\tilde{x}_0)}{4!} \otimes (\tilde{x} \ominus \tilde{x}_0)^4 \oplus \frac{\tilde{f}^5(\tilde{x}_0)}{5!} \otimes (\tilde{x} \ominus \tilde{x}_0)^5 \oplus \ldots . \frac{\tilde{f}^n(\tilde{x}_0)}{n!} \otimes (\tilde{x} \ominus \tilde{x}_0)^n \oplus \ldots$$

Thus, the series approximation of $\tilde{f}$, can be written as,

$$\tilde{f}(\tilde{x}) = \tilde{f}(\tilde{x}_0) \oplus \dot{\tilde{f}}(\tilde{x}_0) \otimes (\tilde{x} \ominus \tilde{x}_0) \oplus \frac{\ddot{\tilde{f}}(\tilde{x}_0)}{2!} \otimes (\tilde{x} \ominus \tilde{x}_0)^2 \oplus \frac{\dddot{\tilde{f}}(\tilde{x}_0)}{3!} \otimes (\tilde{x} \ominus \tilde{x}_0)^3 \oplus \frac{\tilde{f}^4(\tilde{x}_0)}{4!} \otimes (\tilde{x} \ominus \tilde{x}_0)^4 \oplus \frac{\tilde{f}^5(\tilde{x}_0)}{5!} \otimes (\tilde{x} \ominus \tilde{x}_0)^5 \oplus \ldots\ldots\ldots . \frac{\tilde{f}^n(\tilde{x}_0)}{n!} \otimes (\tilde{x} \ominus \tilde{x}_0)^n \oplus \ldots$$

Also, this expansion is true for all *x* belonging to the radius of convergence.

## 3.6 Illustrative examples for Fuzzy Taylor series expansions of fuzzy functions

### 3.6.1 Example:

Fuzzy Taylor series for $\tilde{\boldsymbol{f}}(\tilde{\boldsymbol{x}}) = \boldsymbol{e}^{\tilde{\boldsymbol{x}}}$ about $\tilde{\boldsymbol{x}} = \tilde{\boldsymbol{0}}$ , where $\tilde{\boldsymbol{0}} = (-\boldsymbol{1}, \boldsymbol{0}, \boldsymbol{1})$.
Now for all $\boldsymbol{\alpha}$, ${}^{\boldsymbol{\alpha}}\tilde{\boldsymbol{x}} = {}^{\boldsymbol{\alpha}}\tilde{\boldsymbol{0}} = (\boldsymbol{\alpha} - \boldsymbol{1}, \boldsymbol{1} - \boldsymbol{\alpha})$.

$$\dot{\tilde{f}}(\tilde{x}) = e^{\tilde{x}}, \dot{\tilde{f}}(\tilde{0}) = e^{\tilde{0}} = e^{(\alpha-1,1-\alpha)} = [min(e^{\alpha-1}, e^{1-\alpha}), max(e^{\alpha-1}, e^{1-\alpha})] = [e^{1-\alpha}, e^{\alpha-1}]$$

$$\ddot{\tilde{f}}(\tilde{x}) = e^{\tilde{x}}, \ddot{\tilde{f}}(\tilde{0}) = e^{\tilde{0}} = [e^{1-\alpha}, e^{\alpha-1}]$$

:

:

$$\widetilde{f^n}(\tilde{x}) = e^{\tilde{x}}, \widetilde{f^n}(\tilde{0}) = e^{\tilde{0}} = [e^{1-\alpha}, e^{\alpha-1}]$$

By using these values in fuzzy Taylor's series as in 3.7, we have,

$$e^{\tilde{x}} = e^{\tilde{0}} \oplus e^{\tilde{0}} \otimes (\tilde{x} \ominus \tilde{0}) \oplus \frac{e^{\tilde{0}}}{2} \otimes (\tilde{x} \ominus \tilde{0})^2 \oplus \frac{e^{\tilde{0}}}{3!} \otimes (\tilde{x} \ominus \tilde{0})^3 \oplus \frac{e^{\tilde{0}}}{4!} \otimes (\tilde{x} \ominus \tilde{0})^4 \oplus \frac{e^{\tilde{0}}}{5} \otimes (\tilde{x} \ominus \tilde{0})^5 \oplus \ldots . \frac{e^{\tilde{0}}}{n!} \otimes (\tilde{x} \ominus \tilde{0})^n \oplus \ldots$$

$$e^{\tilde{x}} = e^{\tilde{0}} \otimes [1 \oplus (\tilde{x} \ominus \tilde{0}) \oplus \frac{(\tilde{x}\ominus\tilde{0})^2}{2!} \oplus \frac{(\tilde{x}\ominus\tilde{0})^3}{3!} \oplus \frac{(\tilde{x}\ominus\tilde{0})^4}{4!} \oplus \frac{(\tilde{x}\ominus\tilde{0})^5}{5!} \oplus \ldots \ldots \oplus \frac{(\tilde{x}\ominus\tilde{0})^n}{n!} \oplus \ldots$$

After taking alpha – cut,

$$e^{\tilde{x}} = [e^{1-\alpha}, e^{\alpha-1}][1 \oplus (\tilde{x} \ominus \tilde{0}) \oplus \frac{(\tilde{x}\ominus\tilde{0})^2}{2!} \oplus \frac{(\tilde{x}\ominus\tilde{0})^3}{3!} \oplus \frac{(\tilde{x}\ominus\tilde{0})^4}{4!} \oplus \ldots \ldots \oplus \frac{(\tilde{x}\ominus\tilde{0})^n}{n!} \oplus \ldots]$$

The radius of convergence for $e^{\tilde{x}}$, as defined in 3.4.1,

$$\tilde{R} = \lim_{n\to\infty} \left|\frac{\tilde{a}_n}{\tilde{a}_{n+1}}\right|$$

$$\tilde{R} = \lim_{n\to\infty} \left|\frac{\frac{(\tilde{x}\ominus\tilde{0})^n}{n!}}{\frac{(\tilde{x}\ominus\tilde{0})^{n+1}}{n+1!}}\right|$$

$$\tilde{R} = \lim_{n\to\infty} \left|\frac{(n+1)}{(\tilde{x}\ominus\tilde{0})}\right|$$

So as $n$ tends to infinity, radius of convergence becomes infinite. Thus $e^{\tilde{x}}$ is convergent everywhere.

At core, $\alpha = 1$,

$$e^x = 1 + x + \frac{x^2}{2!} + \frac{x^3}{3!} + \frac{x^4}{4!} + \cdots$$

That is same as Taylor's expansion for real valued function.

### 3.6.2 Example:

Fuzzy Taylor series for $\tilde{f}(\tilde{x}) = sin\tilde{x}$ about $\tilde{x} = \tilde{0}$.

$$\dot{\tilde{f}}(\tilde{x}) = cos(\tilde{x}), \quad \dot{\tilde{f}}(\tilde{0}) = cos\tilde{0}$$
$$\ddot{\tilde{f}}(\tilde{x}) = -sin\tilde{x}, \quad \ddot{\tilde{f}}(\tilde{0}) = -sin\tilde{0}$$
$$\dddot{\tilde{f}}(\tilde{x}) = -cos\tilde{x}, \quad \dddot{\tilde{f}}(\tilde{0}) = -cos\tilde{0}$$
$$\tilde{f}^4(\tilde{x}) = sin\tilde{x}, \quad \tilde{f}^4(\tilde{x}) = sin\tilde{0}$$

:

:

By using these values in fuzzy Taylor's series as in 3.7, we have,

$$sin\tilde{x} = \tilde{f}(\tilde{0}) \oplus \dot{\tilde{f}}(\tilde{0}) \otimes (\tilde{x} \ominus \tilde{0}) \oplus \frac{\ddot{\tilde{f}}(\tilde{0})}{2!} \otimes (\tilde{x} \ominus \tilde{0})^2 \oplus \frac{\dddot{\tilde{f}}(\tilde{0})}{3!} \otimes (\tilde{x} \ominus \tilde{0})^3 \oplus \frac{\tilde{f}^4(\tilde{0})}{4!} \otimes (\tilde{x} \ominus \tilde{0})^4 \oplus \frac{\tilde{f}^5(\tilde{0})}{5!} \otimes (\tilde{x} \ominus \tilde{0})^5 \oplus \ldots \ldots \ldots \ldots \ldots \ldots . \frac{\tilde{f}^n(\tilde{0})}{n!} \otimes (\tilde{x} \ominus \tilde{0})^n \oplus \ldots$$

$$sin\tilde{x} = sin\tilde{0} \oplus cos\tilde{0} \otimes (\tilde{x} \ominus \tilde{0}) \ominus \frac{\sin(\tilde{0})}{2} \otimes (\tilde{x} \ominus \tilde{0})^2 \ominus \frac{\cos(\tilde{0})}{3!} \otimes (\tilde{x} \ominus \tilde{0})^3 + \frac{\sin(\tilde{0})}{4!} \otimes (\tilde{x} \ominus \tilde{0})^4 \oplus \frac{\cos(\tilde{0})}{5!} \otimes (\tilde{x} \ominus \tilde{0})^5 \oplus \ldots$$

The radius of convergence for $sin\tilde{x}$, as defined in 3.4.1,

$$\tilde{R} = \lim_{n\to\infty} \left|\frac{\tilde{a}_n}{\tilde{a}_{n+1}}\right|$$

$$\tilde{R} = \lim_{n\to\infty} \left|\frac{\frac{(\tilde{x}\ominus\tilde{0})^n}{n!}\sin(\tilde{x}+\frac{n\pi}{2})}{\frac{(\tilde{x}\ominus\tilde{0})^{n+1}}{n+1!}\sin(\tilde{x}+\frac{(n+1)\pi}{2})}\right|$$

$$\tilde{R} = \lim_{n\to\infty} \left|\frac{(n+1)}{(\tilde{x}\ominus\tilde{0})}\frac{\sin(\tilde{x}+\frac{n\pi}{2})}{\sin(\tilde{x}+\frac{(n+1)\pi}{2})}\right|$$

So as $n$ tends to infinity, this term $\frac{\sin(\tilde{x}+\frac{n\pi}{2})}{\sin(\tilde{x}+\frac{(n+1)\pi}{2})}$ remains finite and radius of convergence becomes infinite. Thus $sin\tilde{x}$ is convergent everywhere.

At core, $\alpha = 1$,

$$sin = x - \frac{x^3}{3!} + \frac{x^5}{5!} - \frac{x^7}{7!} \ldots$$

**3.6.3 Example:**

Fuzzy Taylor series for $\tilde{f}(\tilde{x}) = cos\tilde{x}$ about $\tilde{x} = \tilde{0}$.

$$\dot{\tilde{f}}(\tilde{x}) = -sin\tilde{x}, \quad \dot{\tilde{f}}(\tilde{0}) = -sin\tilde{0}$$
$$\ddot{\tilde{f}}(\tilde{x}) = -cos\tilde{x}, \quad \ddot{\tilde{f}}(\tilde{0}) = -cos\tilde{0}$$
$$\dddot{\tilde{f}}(\tilde{x}) = sin\tilde{x}, \quad \dddot{\tilde{f}}(\tilde{0}) = sin\tilde{0}$$
$$\tilde{f}^4(\tilde{x}) = cos\tilde{x}, \quad \tilde{f}^4(\tilde{x}) = cos\tilde{0}$$
$$\vdots$$

By using these values in fuzzy Taylor's series as in 3.7, we have,

$$cos\tilde{x} = cos\tilde{0} \ominus sin\tilde{0} \otimes (\tilde{x} \ominus \tilde{0}) \ominus \frac{\cos\tilde{0}}{2!} \otimes (\tilde{x} \ominus \tilde{0})^2 \oplus \frac{sin\tilde{0}}{3!} \otimes (\tilde{x} \ominus \tilde{0})^3 \oplus \frac{cos\tilde{0}}{4!} \otimes (\tilde{x} \ominus \tilde{0})^4 \ominus \frac{\sin\tilde{0}}{5!} \otimes (\tilde{x} \ominus \tilde{0})^5 \oplus \ldots$$

The radius of convergence for $sin\tilde{x}$, as defined in 3.4.1,

$$\tilde{R} = \lim_{n\to\infty} \left|\frac{\tilde{a}_n}{\tilde{a}_{n+1}}\right|$$

$$\tilde{R} = \lim_{n\to\infty} \left|\frac{\frac{(\tilde{x}\ominus\tilde{0})^n}{n!}\cos(\tilde{x}+\frac{n\pi}{2})}{\frac{(\tilde{x}\ominus\tilde{0})^{n+1}}{n+1!}\cos(\tilde{x}+\frac{(n+1)\pi}{2})}\right|$$

$$\tilde{R} = \lim_{n\to\infty} \left|\frac{(n+1)}{(\tilde{x}\ominus\tilde{0})}\frac{\cos(\tilde{x}+\frac{n\pi}{2})}{\cos(\tilde{x}+\frac{(n+1)\pi}{2})}\right|$$

So as $n$ tends to infinity, this term $\frac{\cos(\tilde{x}+\frac{n\pi}{2})}{\cos(\tilde{x}+\frac{(n+1)\pi}{2})}$ remains finite and the radius of convergence becomes infinite. Thus, $cos\tilde{x}$ is convergent everywhere.

At core, $\alpha = 1$,

$$cosx = 1 - \frac{x^2}{2!} + \frac{x^4}{4!} - \frac{x^6}{6!} \dots$$

In the following we solve the example of fully fuzzy differential equation using the fuzzy Taylor's series method.

## 4 Numerical Illustrative

$$\dot{y} = \tilde{x} \otimes \tilde{x} \oplus \tilde{y} \otimes \tilde{y};\ \tilde{y}(\tilde{1}) = \widetilde{2.3} \tag{10}$$

with $\tilde{1} = (0.7,1,1.2)$, $\widetilde{2.3} = (2.1,2.3,2.5)$, $\widetilde{0.1} = (0.07, 0.1, 0.12)$

We solve the above equation for the value $\tilde{y}(\tilde{1} \oplus \widetilde{0.1})$, taking $\tilde{h} = \widetilde{0.1} = (0.07, 0.1, 0.12)$ by proposed fuzzy Taylor's series upto degree 2 as in, Section 3.7. We know,

$$\tilde{y}(\tilde{1} \oplus \widetilde{0.1}) = \tilde{y}(\tilde{1}) \oplus \widetilde{0.1} \otimes \dot{\tilde{y}}(\tilde{1}) \oplus \frac{(\widetilde{0.1})^2}{2!} \otimes \ddot{\tilde{y}}(\tilde{1}) \oplus \dots \tag{11}$$

Now,

$$\tilde{y}(\tilde{1}) = \widetilde{2.3},$$

$$\dot{\tilde{y}}(\tilde{1}) = \tilde{x}^2{}_0 \oplus \tilde{y}^2{}_0 = \tilde{1}^2 \oplus \widetilde{2.3}^2,$$

$$\ddot{\tilde{y}}(\tilde{1}) = 2 \times \left(\tilde{1} \oplus \widetilde{2.3} \otimes (\tilde{1}^2 \oplus \widetilde{2.3}^2)\right),$$

Putting above values in equation (11),

$$\tilde{y}(\tilde{1} \oplus \widetilde{0.1}) = \widetilde{2.3} \oplus \widetilde{0.1} \otimes (\tilde{1}^2 \oplus \widetilde{2.3}^2) \oplus (\widetilde{0.1})^2 \otimes (\tilde{1} \oplus \widetilde{2.3} \otimes (\tilde{1}^2 \oplus \widetilde{2.3}^2)) \oplus \dots \tag{12}$$

Putting the triangular number representation of the fuzzy number involved in (12), we get

$$\tilde{y}(\widetilde{1.01})$$
$$= (2.1, 2.3, 2.5) \oplus (0.07, 0.1, 0.12) \otimes (4.9, 6.29, 7.69) \oplus (0.0049, 0.01, 0.0144) \otimes (10.99, 15.467, 20.425)$$
$$= (2.1, 2.3, 2.5) \oplus (0.343, 0.629, 0.9228) \oplus (0.0538, 0.15467, 0.29412)$$
$$= (2.49, 3.08, 3.71)$$

To visualize the solution in depth, we take the parametric form of the terms in equation (12),

$$^{\alpha}\tilde{y}(^{\alpha}\tilde{1} \oplus {}^{\alpha}\widetilde{0.1})$$
$$= {}^{\alpha}\widetilde{2.3} \oplus {}^{\alpha}\widetilde{0.1} \otimes ({}^{\alpha}\tilde{1}^2 \oplus {}^{\alpha}\widetilde{2.3}^2) \oplus \widetilde{0.1}^2 \otimes \left({}^{\alpha}\tilde{1} \oplus {}^{\alpha}\widetilde{2.3} \otimes ({}^{\alpha}\tilde{1}^2 \oplus {}^{\alpha}\widetilde{2.3}^2)\right) \tag{13}$$

The parametric form are as follows,

$$^{\alpha}\tilde{y} = [\underline{y}, \overline{y}],$$

$$^{\alpha}\widetilde{2.3} = [2.1 + 0.2\alpha, 2.5 - 0.2\alpha],$$

$${}^{\alpha}\tilde{1} = [0.7 + 0.3\alpha, 1.2 - 0.2\alpha],$$

$${}^{\alpha}\widetilde{0.1} = [0.07 + 0.03\alpha, 0.12 - 0.02\alpha].$$

Putting the parametric values in (13) we get,

$$[\underline{y}({}^{\alpha}\tilde{1} \oplus {}^{\alpha}\widetilde{0.1}), \overline{y}({}^{\alpha}\tilde{1} \oplus {}^{\alpha}\widetilde{0.1})]$$

$$= [2.1 + 0.2\alpha, 2.5 - 0.2\alpha]$$

$$+[3.9 \times 10^{-3}\alpha^3 + 0.0469\alpha^2 + 0.2352\alpha + 0.343, -1.6 \times 10^{-3}\alpha^3 + 0.0392\alpha^2 - 0.3314\alpha + 0.9228]$$

$$+[(0.07 + 0.03\alpha)^2,\ (0.12 - 0.02\alpha)^2]\,[0.026\alpha^3 + 0.525\alpha^2 + 3.926\alpha + 10.99, -0.016\alpha^3 + 0.496\alpha^2 - 2.478\alpha + 20.425]\ldots$$

$$\therefore \left[\underline{y}([0.7 + 0.3\alpha, 1.2 - 0.2\alpha] + [0.07 + 0.03\alpha, 0.12 - 0.02\alpha]), \overline{y}([0.7 + 0.3\alpha, 1.2 - 0.2\alpha] + [0.07 + 0.03\alpha, 0.12 - 0.02\alpha])\right]$$

$$= [2.1 + 0.2\alpha, 2.5 - 0.2\alpha]$$

$$+[3.9 \times 10^{-3}\alpha^3 + 0.0469\alpha^2 + 0.2352\alpha + 0.343, -1.6 \times 10^{-3}\alpha^3 + 0.0392\alpha^2 - 0.3314\alpha + 0.9228]$$

$$+[(0.07 + 0.03\alpha)^2,\ (0.12 - 0.02\alpha)^2]\,[0.026\alpha^3 + 0.525\alpha^2 + 3.926\alpha + 10.99, -0.016\alpha^3 + 0.496\alpha^2 - 2.478\alpha + 20.425].$$

$$\therefore \left[\underline{y}[0.77 + 0.33\alpha, 1.32 - 0.22\alpha]\,, \overline{y}[0.77 + 0.33\alpha, 1.32 - 0.22\alpha]\right]$$

$$= [2.1 + 0.2\alpha, 2.5 - 0.2\alpha]$$

$$+[3.9 \times 10^{-3}\alpha^3 + 0.0469\alpha^2 + 0.2352\alpha + 0.343, -1.6 \times 10^{-3}\alpha^3 + 0.0392\alpha^2 - 0.3314\alpha + 0.9228]$$

$$+[(0.07 + 0.03\alpha)^2,\ (0.12 - 0.02\alpha)^2]\,[0.026\alpha^3 + 0.525\alpha^2 + 3.926\alpha + 10.99, -0.016\alpha^3 + 0.496\alpha^2 - 2.478\alpha + 20.425].$$

Now, we can write left side of the above equation as $\left[\underline{\underline{y}}, \overline{\overline{y}}\right]$ defined in Section 3.3.

Where, $\underline{\underline{y}} = \underline{y}\,[0.77 + 0.33\alpha, 1.32 - 0.22\alpha]$ and $\overline{\overline{y}} = \overline{y}\,[0.77 + 0.33\alpha, 1.32 - 0.22\alpha]$.

So, the equation becomes,

$$\left[\underline{\underline{y}}, \overline{\overline{y}}\right] = [2.1 + 0.2\alpha, 2.5 - 0.2\alpha]$$

$$+[3.9 \times 10^{-3}\alpha^3 + 0.0469\alpha^2 + 0.2352\alpha + 0.343, -1.6 \times 10^{-3}\alpha^3 + 0.0392\alpha^2 - 0.3314\alpha + 0.9228]$$

$$+[(0.07 + 0.03\alpha)^2,\ (0.12 - 0.02\alpha)^2]\,[0.026\alpha^3 + 0.525\alpha^2 + 3.926\alpha + 10.99, -0.016\alpha^3 + 0.496\alpha^2 - 2.478\alpha + 20.425].$$

At, $\alpha = 0$ , $\left[\underline{\underline{y}}, \overline{\overline{y}}\right] = [2.49, 3.71]$ and at $\alpha = 1$, $y(1.1) = 3.08$.

Thus, $\tilde{y}(\widetilde{1.1}) = (2.49\,, 3.08, 3.71) = \widetilde{3.08}$.

Also, it can be shown that the crisp solution matches with ${}^{1}\tilde{y}$, considering the constants involved as crisp reals.

## 5 Conclusion

In this article, we have proposed fuzzy Taylor's series under Modified Hukuhara derivative. We have proposed and proved the technique to obtain the solution in a fully fuzzy environment along with its convergence. The advantage of this fuzzy Taylor series expansion is that we can directly solve fuzzy differential equations in a fuzzy environment without converting them into their crisp counterpart. For future scope, one can develop other techniques completely in a fuzzy environment.